\documentclass[10pt]{amsart}
\usepackage{color}

\definecolor{c20}{rgb}{0.,0.7,0.}
\definecolor{c30}{rgb}{0.,0.,1.}
\definecolor{c40}{rgb}{0.3,0.3,0.9}
\definecolor{c50}{rgb}{1,0,0}
\definecolor{c60}{rgb}{1,0.9,0.1}

\def\gE#1{\textcolor{c40}{#1}}
\def\gE#1{#1}

\def\ze#1{\textcolor{c20}{#1}}
\def\ze#1{#1}
\newcommand{\tb}[1]{{\textcolor{blue}{#1}}}
\def\tb#1{#1}

\def\tb#1{#1}

\def\cLP#1{\textcolor{c50}{#1}}
\def\cLP#1{#1}

\def\cP#1{\textcolor{c50}{#1}}
\def\cP#1{#1}

\newcommand{\kb}[1]{\boldsymbol{#1}}
\newcommand{\vk}[1]{\kb{#1}}

\newcommand{\abs}[1]{\left\lvert #1 \right\rvert}

\newcommand{\E}[1]{\mathbb{E}\left\{#1\right\}}

\newcommand{\pk}[1]{\mathbb{P} \left\{ #1 \right \} }

\newcommand{\R}{\mathbb{R}}

\newcommand{\N}{\mathbb{N}}
\newcommand{\inr}{\in \R}
\newcommand{\inn}{\in \N}
\newcommand{\ldot}{,\ldots,}

\newcommand{\limit}[1]{\lim_{#1 \to   \infty}}

\newcommand{\BQN}{\begin{eqnarray}}
\newcommand{\EQN}{\end{eqnarray}}
\newcommand{\BQNY}{\begin{eqnarray*}}
\newcommand{\EQNY}{\end{eqnarray*}}

\newcommand{\BS}{\begin{sat}}
\newcommand{\ES}{\end{sat}}
\newcommand{\BT}{\begin{theo}}
\newcommand{\ET}{\end{theo}}
\newcommand{\BK}{\begin{korr}}
\newcommand{\EK}{\end{korr}}

\newcommand{\BD}{\begin{de}}
\newcommand{\ED}{\end{de}}

\newcommand{\BIT}{\begin{itemize}}
\newcommand{\EIT}{\end{itemize}}
\newcommand{\BDI}{\begin{description}}
\newcommand{\EDI}{\end{description}}

\newcommand{\BRM}{\begin{remarks}}
\newcommand{\ERM}{\end{remarks}}

\newcommand{\BEL}{\begin{lem}}
\newcommand{\EEL}{\end{lem}}

\newtheorem{theo}{Theorem}[section]
\newtheorem{sat}[theo]{Proposition}
\newtheorem{de}[theo]{Definition}
\newtheorem{lem}[theo]{Lemma}

\newtheorem{korr}[theo]{Corollary}

\newtheorem{remarks}[theo]{Remarks}

\newcommand{\nelem}[1]{{Lemma \ref{#1}}}

\newcommand{\netheo}[1]{{Theorem \ref{#1}}}

\newcommand{\prooftheo}[1]{ \textsc{\bf Proof of Theorem} \ref{#1}:}

\newcommand{\prooflem}[1]{\textsc{\bf Proof of Lemma} \ref{#1}:}

\newcommand{\COM}[1]{}

\newcommand{\QED}{\hfill $\Box$}

\def\rw{\rightarrow}

\def\IF{\infty}

\def\LT{\left}
\def\RT{\right}
\def\H{\mathcal{H}}

\def\ooo{(1+o(1))}

\def\rw{\rightarrow}

\def\H{\mathcal{H}}

\def\almin{\alpha_{\min}}
\newcommand{\todis}{\stackrel{d}{\to}}

\begin{document}

\title{\ze{On the Probability of Conjunctions} of Stationary Gaussian Processes}

\author{Krzysztof D\c{e}bicki}
\address{Krzysztof D\c{e}bicki, Mathematical Institute, University of Wroc\l aw, pl. Grunwaldzki 2/4, 50-384 Wroc\l aw, Poland}
\email{Krzysztof.Debicki@math.uni.wroc.pl}

\author{Enkelejd  Hashorva}
\address{Enkelejd Hashorva, Department of Actuarial Science, 
University of Lausanne,\\
UNIL-Dorigny, 1015 Lausanne, Switzerland
}
\email{Enkelejd.Hashorva@unil.ch}

\author{Lanpeng Ji}
\address{Lanpeng Ji, Department of Actuarial Science, 
University of Lausanne\\
UNIL-Dorigny, 1015 Lausanne, Switzerland
}
\email{Lanpeng.Ji@unil.ch}

\author{Kamil Tabi\'{s}}
\address{Kamil Tabi\'{s}, Department of Actuarial Science, 
University of Lausanne, UNIL-Dorigny, 1015 Lausanne, Switzerlandand
and
Mathematical Institute, University of Wroc\l aw, pl. Grunwaldzki 2/4, 50-384 Wroc\l aw, Poland
}
\email{Kamil.Tabis@unil.ch}

\bigskip

\date{\today}
 \maketitle

{\bf Abstract:} Let $\{X_i(t),t\ge0\}, 1\le i\le n$ be independent
centered stationary Gaussian processes with unit variance and
almost surely continuous sample paths. For given positive
constants $u,T$, define the set of conjunctions $ C_{[0,T],u}:=
\{t\in [0,T]: \min_{1 \le i \le n} X_i(t) \ge u\}.$
\tb{Motivated by
some applications in brain mapping and digital communication systems,
we obtain exact asymptotic expansion of
$
\pk{C_{[0,T],u} \not=\phi},
$
as $u\to\infty$.
Moreover, we establish the Berman sojourn limit theorem for the random process
$\{\min_{1 \le i \le n} X_i(t), t\ge0\}$
and derive the tail asymptotics of the supremum of \ze{each} order statistics process.
}\\



\section{Introduction \& Main Result}
Let  $X_i(\vk{t})$ model the value of an image $i$ at location $\vk{t}\inr^d$,  $1\le i\le n$.
 For a given positive threshold $u$ and a given scan set $\mathcal{T} \subset \R^d$,
the set of conjunctions $C_{\mathcal{T},u}$ is defined by
$$ C_{\mathcal{T},u}:= \{\vk{t}\in \mathcal{T}: \min_{1 \le i \le n} X_i(\vk{t}) \ge u\}$$
see the seminal contribution \cite{MR1747100}. As mentioned in the aforementioned paper, of interest is the \ze{calculation of the} probability that the set of conjunctions $C_{\mathcal{T},u}$ is not empty, i.e., 
$$ p_{\mathcal{T},u}:=\pk{C_{\mathcal{T},u} \not=\phi}= \pk{\sup_{\vk{t}\in \mathcal{T}} \min_{1 \le i \le n} X_i(\vk{t}) \ge u}.$$
Typically, in applications such as the
analysis of functional magnetic resonance imaging (fMRI) data,  $X_i$'s are assumed to be real-valued Gaussian random fields.
Approximations of $p_{\mathcal{T},u}$ are discussed for smooth  Gaussian random fields in \cite{MR1747100, MR2775212, ChengXiao13};
results for non-Gaussian random fields can be found in \cite{MR2654766}.\\
In this paper, we shall consider the case $d=1,\mathcal{T}:=[0,T]$, with $T>0$, and that $X_i$'s are independent centered stationary Gaussian processes with unit  variance and correlation functions $r_i(\cdot), 1\le i\le n$ that satisfy 
\BQN\label{eq:rr}
 r_i(t)= 1 - C_i \abs{t}^{\alpha_i}+ o(\abs{t}^{\alpha_i}), \quad t\to 0, \quad \quad r_i(t)< 1, \quad \forall t\in (0,T]
\EQN
for some positive constants $\alpha_i \in (0,2]$ and $C_i, 1\le i\le n$. Further, we assume that $X_i$'s have almost surely continuous sample paths.
Since the calculation of $p_{\mathcal{T},u}$ is not possible in general, we shall investigate below the exact asymptotic behavior of
$p_{\mathcal{T},u}$ as $u\to \IF$.
Although $\{\min_{1 \le i \le n} X_i(t), t\ge0\}$ is not a Gaussian process when
$n\ge 2$, as shown in \cite{MR1747100},
it happens that it is possible to \tb{adapt} techniques used in the theory of Gaussian processes and random fields to this class of processes.
Motivated by a recent paper of Albin and Choi \cite{AlbinC} and the extremal theory for stationary Gaussian processes developed mainly
by Berman and Albin (see \cite{Berman82,Berman92,MR1043939,Albin2003}), we shall derive an asymptotic expansion for $p_{\mathcal{T},u}$ as $u\to\IF$, by following the ideas of \cite{AlbinC}.

For the formulation of our main result we need to introduce some notation. Let
\tb{ $\{B_{\alpha_i}(t),t\ge 0\}, 1\le i\le n$  be mutually independent standard fractional Brownian motions with Hurst
indexes $\alpha_i/2\in(0,1]$, $1\le i\le n$, respectively,
i.e., $B_{\alpha_i}$ is a centered Gaussian process with continuous sample paths and
covariance function
}
$$Cov(B_{\alpha_i}(t), B_{\alpha_i}(s))= \frac{1}{2} \Bigl(t^{\alpha_i}+ s^{\alpha_i} - \abs{t-s}^{\alpha_i}\Bigr), \quad t,s>0,\ \ \ 1\le i\le n.$$
 Next define 
\BQN\label{mathcZ}
\mathcal{Z}(t):= \min_{1 \le i \le n} \Bigl( \Bigl(\sqrt{2} B_{\alpha_i}(C_i^{1/\alpha_i}t)
- C_i t^{\alpha_i}\Bigr) \vk{1}(\alpha_i= \almin)+ E_i\Bigr),\ \ t\ge0, \  \  \alpha_{\min}:= \min_{1 \le i \le n}  \alpha_i,
\EQN
where $\vk{1}(\cdot)$ denotes the indicator function, and ${E_i}$'s are
\tb{mutually}
independent unit exponential random variables being further independent of ${B_{\alpha_i}}$'s.
\tb{Finally, let
$\mathcal{H}_{\alpha_1 \ldot \alpha_n}(C_1\ldot C_n) \in (0,\IF)$
denote
{\it generalized Pickands constant}, determined by
\BQN\label{HAAN}
\mathcal{H}_{\alpha_1 \ldot \alpha_n}(C_1\ldot C_n)= \lim_{a \downarrow 0} \frac{1}{a}\pk{ \max_{k\ge 1} \mathcal{Z}(ak)\le 0}.
\EQN
The following theorem constitutes our principle result.
}
\BT\label{ThmB} Let $\{X_i(t), t\ge0\}$, $1\le i\le n$ be \tb{mutually}
independent centered stationary Gaussian processes with 
unit variance and correlation functions satisfying \eqref{eq:rr}.
Then, for any $T>0$
\BQN \label{eqMain}
\quad \pk{\sup_{t \in [0,T]} \min_{1 \le i \le n} X_i(t)> u}=\mathcal{H}_{\alpha_1 \ldot \alpha_n}(C_1\ldot C_n)\ T  u^{\frac{2}{\almin}} \frac{\exp(- n u^2/2)}{ (2 \pi)^{n/2} u^n}  
(1+ o(1)), \quad u\to \IF,
\EQN
\ze{where $\mathcal{H}_{\alpha_1 \ldot \alpha_n}(C_1\ldot C_n) \in (0,\IF)$ is defined in \eqref{HAAN}.}
\ET

The organization of the paper. Section 2 presents brief discussions and shows the validity of
the Berman sojourn limit theorem for the random process $\{\min_{1 \le i \le n} X_i(t), t\ge0\}$.
\tb{Additionally, utilizing the fact that
minimum is a particular case of the order statistics, in Theorem \ref{ThmB2} we get a counterpart of
Theorem \ref{ThmB} for order statistics processes.
The case of non-standard stationary Gaussian processes is treated in \netheo{ThmNS}.}
Section 3 contains all the proofs.

\section{Discussions \& Extensions}
In his seminal contribution
\cite{PicandsA} J. Pickands III established the exact asymptotic tail behavior of the supremum of the stationary Gaussian process $\{X_1(t),t\in[0,T]\}$ under the condition \eqref{eq:rr}, using a double-sum method. The first crucial step to that result
is the celebrated Pickands lemma which states that, for any positive constant $S$
\BQNY 
\pk{\sup_{t\in[0,u^{-\cLP{\frac{2}{\alpha_1}}}S]}X_1(t)>u}=\mathcal{H}_{\alpha_1}\LT[0, C_1^{\frac{1}{\alpha_1}}S\RT]\Psi(u)  \ooo, \ \ \ \ \ u \rw \IF,
\EQNY
where $\Psi(\cdot)$ is the survival function of an  $N(0,1)$ random variable and
\BQNY
\mathcal{H}_{\alpha_1}[0, S]=\E{ \exp\biggl(\sup_{t\in[0,S]}\Bigl(\sqrt{2}B_{\alpha_1} (t)-t^{\alpha_1}\Bigr)\biggr)}\in(0,\IF).
 \EQNY
Recall that $\Psi(u)= \exp(- u^2/2)/ \sqrt{2 \pi u^2}(1+o(1))$ as $u\to \IF$.\\
An application of Pickands lemma, together with the double-sum method,
\tb{yields}
(see, e.g., \cite{PicandsA, Pit72, Pit96})
\BQN \label{pic}
\pk{\sup_{t\in[0,T]}X_1(t)>u}=T C_1^{\frac{1}{\alpha_1}} \mathcal{H}_{\alpha_1} u^{\frac{2}{\alpha_1}}\Psi(u)\ooo, \quad u \rw \IF,
\EQN
where $\mathcal{H}_{\alpha_1}\in (0,\IF)$  is the Pickands constant, defined by
\BQNY
\mathcal{H}_{\alpha_1}=\lim_{S\to\IF}\frac{1}{S}\mathcal{H}_{\alpha_1}[0, S].
\EQNY
\tb{
We refer to the recent contribution \cite{DikerY}, where alternative
representations of Pickands constant are derived; see also \cite{debicki2002ruin,Faletal2010}
and the references therein for properties and generalizations of Pickands constant.
}\\
The constant $\mathcal{H}_{\alpha_1 \ldot \alpha_n}(C_1\ldot C_n)$, appearing in \eqref{HAAN},
is more complicated than $\mathcal{H}_{\alpha_1}$.
A \tb{simple} lower bound for $\mathcal{H}_{\alpha_1 \ldot \alpha_n}(C_1\ldot C_n)$ can be found as follows:
 \BQN\label{lowerB}
 \mathcal{H}_{\alpha_1 \ldot \alpha_n}(C_1\ldot C_n)&\ge & \max_{1 \le i \le n: \alpha_i=\almin}
 \lim_{a \downarrow 0} \frac{1}{a}\pk{ \max_{k\ge \ze{1}}  \Bigl(\sqrt{2} B_{\alpha_i}(C_i^{1/\alpha_i}a k)
- C_i (a k)^{\alpha_i}\Bigr) + E_i \le 0}\notag\\
&\ge & \max_{1 \le i \le n: \alpha_i=\almin} C_i^{1/\alpha_i}
 \lim_{a \downarrow 0} \frac{1}{C_i^{1/\alpha_i}a }\pk{ \max_{k\ge \ze{1}}(\sqrt{2} B_{\alpha_i}(C_i^{1/\alpha_i}a k)
-  (C_i^{1/\alpha_i}a k)^{\alpha_i})+ E_i\le0}\notag\\
&\ge & \max_{1 \le i \le n: \alpha_i=\almin} C_i^{1/\almin} \mathcal{H}_{\almin}>0,
 \EQN
 where in the last step we used the alternative expression of the Pickands constant given in \cite{AlbinC}.

\tb{\netheo{ThmB} can also} be proved using the double-sum method, extending thus the Pickands lemma and Pickands theorem to include the non-Gaussian process $\{\min_{1 \le i \le n} X_i(t),t\ge 0\}$; due to heavy technical details the proof will be displayed in a forthcoming article.


Finally, we remark that in view of the recent contributions  \cite{MR2733939,Turkman2012}
 it is possible to derive the exact asymptotics of $p_{\mathcal{T},u}$
considering $X_i(\vk{t}),\vk{t}\inr^d$ stationary isotropic Gaussian random fields.

We continue below with \cP{four results, the first one establishes a Berman sojourn limit theorem,
the second one deals with order statistics processes of $X_i$'s, the third one focuses
on \ze{a} time-changed model, and the last one concerns
a generalization of \netheo{ThmB} to non-standard stationary Gaussian $X_i$'s.}
\tb{
\subsection{A  \gE{Berman} sojourn limit theorem}
Let, \gE{for $t\ge 0$},
\BQN\label{eq:soj}
L_t(u)= \int_0^t \vk{1}(\min_{ 1\le i\le n} X_i(s)> u) \, ds
 \EQN
be the sojourn time of the process $\{\min_{1\le i\le n} X_i(t), t\ge0\}$ above a level $u>0$ on the time interval $[0,t]$.
}
The next result is the Berman sojourn limit theorem for the process $\{\min_{ 1\le i\le n} X_i(t), t\ge0\}$.

\BT \label{Thm:Berman}
Let $\{X_i(t),t\ge0\}, 1\le i\le n$ be independent centered  stationary Gaussian processes  with unit variance and correlation functions that satisfy  \eqref{eq:rr}, and let $L_t(u)$ be defined as in \eqref{eq:soj} for any positive constants $t,u$. Then we have,
for all $t>0$ small enough, that
$$ \limit{u}  \int_x^\IF  \frac{ \pk{u^\frac{2}{\almin}L_t(u)> y  }}{ u^{\frac{2}{\almin}}\E{ L_t(u)}}\, dy =  B(x) $$
holds at all continuity points $x>0$ of  $B(x)= \pk{\int_0^\IF \vk{1}(\mathcal{Z}(s)> 0)\, ds > x}$.
\ET
\tb{
\subsection{Asymptotics of supremum of order statistics processes}
Let $\{X_{i:n}(t),t\ge0\}, 1\le i\le n$
be the order statistics processes of $\{X_i(t),t\ge0\}, 1\le i\le n$, i.e., we \gE{define}
}
\BQNY 
X_{1:n}(t):=\max_{1 \le i \le n} X_i(t)\ge X_{2:n}(t)\ge \ldots \ge X_{n:n}(t)=\min_{1 \le i \le n} X_i(t),\ \ \ t\ge0.
\EQNY
\tb{
Our next result concerns the exact tail asymptotics of the supremum of the order statistics processes.
We refer, e.g., to \cite{SaKa05} for motivation of study the exit probabilities of the order statistics processes
in electrical engineering.
}

For \tb{clearness of the presentation}, we assume further that $\alpha_1=\ldots=\alpha_n=:\alpha$ and $C_1=\ldots=C_n=1$. Furthermore, define
\BQNY 
\mathcal{H}_{\alpha,j}= \lim_{a \downarrow 0} \frac{1}{a}\pk{ \max_{k\ge 1} \mathcal{Z}_j(ak)\le 0},\ \ 1\le j\le n,
\EQNY
where
$$ \mathcal{Z}_j(t):= \min_{1 \le i \le j} \Bigl(  \sqrt{2} B_{\alpha}^{(i)}(t)
-  t^{\alpha} + E_i\Bigr),\ \ t\ge0, $$
with  ${E_i}$'s being independent unit exponential random variables  which are further independent of mutually independent fractional Brownian motions  $B_{\alpha}^{(i)}$'s.

\BT\label{ThmB2}
Let $\{X_i(t),t\ge0\}, 1\le i\le n$ be independent centered  stationary Gaussian processes
with unit variance and correlation functions that satisfy  \eqref{eq:rr} with
$\alpha_1=\ldots=\alpha_n=:\alpha$
and
$C_1=\ldots=C_n=1$.
\tb{Then, for any $T>0$}
\BQNY 
\pk{\sup_{t \in [0,T]} X_{j:n}(t)> u}= \mathcal{H}_{\alpha, j}\ T \frac{n!}{(n-j)!j!} u^{\frac{2}{\alpha}}  (\Psi(u))^j(1+ o(1)), \quad  1\le j\le n
\EQNY
\tb{
as $u\to \IF$.}
\ET
\tb{
\subsection{Conjunction of time-changed processes}
The technique of Albin and Choi \cite{AlbinC} which we applied in the proof of Theorem \ref{ThmB},
can be utilized also for some other interesting extensions.
}
To illustrate it, we \gE{investigate} the tail asymptotics of supremum of process  $Y_{\vk{\Theta}}(t)= \min_{1 \le i \le n} X_i^*(t), t\ge 0$ on a finite-time interval, say $[0,T]$, where $X_i^*(t)= X_i(\Theta_i t), 1\le i\le n$ are time-changed centered stationary Gaussian processes with $\Theta_i$'s  non-degenerate non-negative bounded random variables being independent of $X_i$'s; \gE{see \cite{AREN,DEJ13} for recent results on the extremes of time-changed Gaussian processes.} Indeed, it follows easily that the result of Lemma 3.1 (see Section 3) holds with limit process
\BQNY
\mathcal{Z}_{\vk{\Theta}}(t)= \min_{1 \le i \le n} \Bigl(  \Bigl(\sqrt{2} B_{\alpha_i}(C_i^{1/\alpha_i} \Theta_i t)
- C_i (\Theta_i t)^{\alpha_i}\Bigr) \vk{1}(\alpha_i= \almin )+ E_i\Bigr), \ \ t\ge0,
\EQNY
where ${B_{\alpha_i}}$'s and ${E_i}$'s are given as before which are further independent of $\Theta_i$'s.
Thus, we have by a similar proof as \netheo{ThmB} that
\BQNY 
\pk{\sup_{t \in [0,T]} Y_{\vk{\Theta}}(t)> u}= \mathcal{H}_{\alpha_1 \ldot \alpha_n}^*(C_1\ldot C_n)\ T u^{\frac{2}{\almin}} (\Psi(u))^n(1+ o(1)) , \quad u\to \IF,
\EQNY
where
$$ \mathcal{H}_{\alpha_1 \ldot \alpha_n}^*(C_1\ldot C_n)= \lim_{a \downarrow 0} \frac{1}{a}\pk{ \max_{k\ge 1} \mathcal{Z}_{\vk{\Theta}}(ak)\le 0} \in (0,\IF). $$

\subsection{Non-standard stationary Gaussian processes}
\cP{Let $\widetilde{ X_i}(t)= X_i(t)/ b_i, t\ge0$ for some $b_i>0, 1\le i\le n$
with $X_i$'s being given as in \netheo{ThmB}.
Clearly,  $\widetilde{ X_i}$'s are  again  centered stationary Gaussian processes.
We have the following result considering the supremum of $\min_{1 \le i \le n} \widetilde{ X_i}(t),t\in [0,T]$.
\BT\label{ThmNS}
Under the assumptions of \netheo{ThmB}, we have, for any $T>0$,
\BQN \label{eq:NS}
\quad \pk{\sup_{t \in [0,T]} \min_{1 \le i \le n} \widetilde{X_i}(t)> u}=\widetilde{\mathcal{H}}_{\alpha_1 \ldot \alpha_n}(C_1\ldot C_n)\ T  u^{\frac{2}{\almin}} \prod_{i=1}^n  \Psi(b_i u)
(1+ o(1)), \quad u\to \IF,
\EQN
where
\BQNY
\widetilde{\mathcal{H}}_{\alpha_1 \ldot \alpha_n}(C_1\ldot C_n)= \lim_{a \downarrow 0} \frac{1}{a}\pk{ \max_{k\ge 1} \widetilde{\mathcal{Z}}(ak)\le 0},
\EQNY
with
\BQNY
\widetilde{\mathcal{Z}}(t):= \min_{1 \le i \le n} \LT( \Bigl(\sqrt{2}  b_i^{-1} B_{\alpha_i}(C_i^{1/\alpha_i}t)
- C_i t^{\alpha_i}\Bigr) \vk{1}(\alpha_i= \almin)+ b_i^{-2} E_i\RT),\ \ t\ge0,
\EQNY
and $B_{\alpha_i}$'s and $E_i$'s being given as in Section 1.
\ET
}

\section{Proofs}
\newcommand{\equaldis}{\stackrel{d}{=}}
\tb{The idea of the proof of Theorem \ref{ThmB} is based on the technique developed by
Albin and Choi \cite{AlbinC}. We begin with
several lemmas for the minimum process $Y(t):=X_{n:n}(t)=\min_{1 \le i \le n} X_i(t)$,
which altogether will be used to show the proof of \netheo{ThmB}.
}
Then we present the proofs of Theorem \ref{Thm:Berman} and \netheo{ThmB2}.

Hereafter we shall use the notation and the assumptions of Introduction and Section 2. For notational simplicity we shall set
below
$$q(u)=u^{-2/\almin}, \quad u>0$$
and shall use the standard notation  $\lfloor\cdot\rfloor$ for the the ceiling function, i.e., $\lfloor x\rfloor$ is the largest integer that is smaller than $x\inr$.

\BEL \label{AlbA}
For any grid of points $0\le t_0 < t_1 < \cdots < t_d< \IF$, $d\inn$, we have
the joint convergence in distribution
\BQNY
\Bigl( nu(Y(q(u) t_1)- u) \ldot nu(Y(q(u)t_d)- u) \Bigr) \Bigl \lvert (Y(0)>u) &\todis & n \Bigl( \mathcal{Z}(t_1) \ldot \mathcal{Z}(t_d)\Bigr)
\EQNY
as $u\to\IF$, where the process $\mathcal{Z}$ is defined as in \eqref{mathcZ}.
\EEL

\prooflem{AlbA} First note that $Y(0)$ has distribution function $G(\cdot)$ in the Gumbel max-domain of attraction with positive scaling function $w(u)=n u$ i.e.,
$$ \limit{u} \frac{1-G( u+  {x}/{w(u)} )}{1-G(u)}=\exp(-x), \quad x\inr. $$
See \cite{Faletal2010, Res1987} for more details on the Gumbel max-domain of attraction. Moreover, it follows from
Lemma 2 in \cite{AlbinC}
that, for any $1\le i\le n$, the
following joint convergence in distribution
\BQNY
\Bigl( X_{iu}(t_1 ) \ldot X_{iu}(t_d) \Bigr) \biggl \lvert (X_{i}(0)>u) \todis  \biggl( \sqrt{2} B_{\alpha_i}(C_i^{1/\alpha_i}t_1)- C_it_1^{\alpha_i}+ E_i \ldot \sqrt{2}B_{\alpha_i}(C_i^{1/\alpha_i}t_d)- C_it_d^{\alpha_i} + E_i \biggr)
\EQNY
holds as $ u\to \IF$, where $X_{iu}(t)= u (X_{i}( u^{-2/\alpha_i} t)- u), t\ge0, u>0$.
Then the claim follows by the independence of  $X_i$'s, $B_{\alpha_i}$'s and $E_i$'s.  \QED

\BEL \label{AlbB}  For any $a>0$ we have
\BQNY
\limit{N} \limit{u} \frac{1}{N \pk{Y(0)> u}} \pk{\max_{k \in \{0 \ldot N\}} Y(a q(u) k) > u} &=&
\pk{ \bigcap_{l=1}^\IF \{\mathcal{Z}(al) \le 0\}}.
\EQNY
\EEL
\prooflem{AlbB} In view of Lemma \ref{AlbA}, the proof follows with the same arguments as that of Lemma 3 in \cite{AlbinC}. \QED


\def\oA{\overline{\alpha}}
\BEL \label{AlbC}  
Let $\oA:=\almin/4$. We have
\BQNY
\lim_{ a \downarrow 0} \limsup_{u\to \IF} \frac{q(u)}{\pk{Y(0)>u}} \pk{
\sup_{t\in [0,T]} Y(t)> u + \frac{a^{\oA}}{u},
\max_{ k\in \{0 \ldot  \lfloor T /(aq(u)) \rfloor\}} Y(aq(u) k) \le u}=0.
\EQNY
\EEL
\prooflem{AlbC} The proof follows by  similar arguments as that of Lemma 4 in  \cite{AlbinC}.  Since the proof of Lemma 4 in \cite{AlbinC} only requires the stationarity and the continuity of the process involved, we obtain,
for all large $u$ and small $a>0$, that
\BQNY
\lefteqn{\pk{\sup_{t\in [0,T]} Y(t)> u + \frac{a^{\oA}}{u},
\max_{ k\in \{0 \ldot  \lfloor T a/q(u) \rfloor\}} Y(aq(u) k)\le u}}\\
 &\le & \frac{2 T}{a q(u)} \sum_{j=1}^\IF 2^{ j}
\pk{ Y(a q(u)  2^{-j}) > u_j, Y(0) \le u_{j-1} },
\EQNY
where $u_j:=u + a^{\oA}(1- 2^{-j \oA})/u >u$ for any $j\ge1$. Further,  for all $u$ large and $a>0$ small and any $1\le i\le n, j\ge 1$, the following inequality
\BQNY
r_i( a q(u)  2^{-j} ) X_i(a q(u)  2^{-j}) - X_i(0) &\ge& a^{\oA}(2^{\oA}-1) 2^{-j \oA-1}/u=:c_{ju}
\EQNY
is implied by the event $\{X_i(a q(u)  2^{-j}) > u_j, X_i(0) \le u_{j-1}\}$.
Thus, in view of  the fact that $r_i( a q(u)  2^{-j} ) X_i(a q(u)  2^{-j}) - X_i(0) $
is independent of
\tb{$X_i(a q(u)  2^{-j})$,}
we conclude that
\BQNY
\lefteqn{\pk{\sup_{t\in [0,T]} Y(t)> u + \frac{a^{\oA}}{u},
\max_{ k\in \{0 \ldot  \lfloor T a/q(u) \rfloor\}} Y(aq(u) k) \le u}}\\
 &\le & \frac{2 T}{a q(u)} \sum_{j=1}^\IF 2^{ j}
\pk{ X_1 (a q(u)  2^{-j}) > u_j, \ldots, X_n (a q(u)  2^{-j}) > u_j , \bigcup_{i=1}^n \{X_i(0) \le u_{j-1} \}}\\
&\tb{\le}& \frac{2 T}{a q(u)} \sum_{j=1}^\IF 2^{ j}
\pk{ X_1 (a q(u)  2^{-j}) > u_j, \ldots, X_n (a q(u)  2^{-j}) > u_j , \bigcup_{i=1}^n \{r_i( a q(u)  2^{-j} ) X_i(a q(u)  2^{-j}) - X_i(0) \ge c_{ju}\}}\\
&\le& \frac{2 T}{a q(u)} \sum_{j=1}^\IF 2^{ j} (\Psi(u))^n
\pk{ \bigcup_{i=1}^n \{r_i( a q(u)  2^{-j} ) X_i(a q(u)  2^{-j}) - X_i(0) \ge c_{ju}\}}\\
&\le& \frac{2 T}{a q(u)} \sum_{j=1}^\IF 2^{ j} (\Psi(u))^n  \sum_{i=1}^n \Psi\biggl( \frac{c_{ju}} { \sqrt{1 - r_i( a q(u)  2^{-j} )^2 }} \biggr)
\EQNY
for all $u$ large and $a>0$ small, hence the claim follows. \QED

\prooftheo{ThmB}  The proof
\tb{follows by a similar idea as used in the proof of Theorem 1 in \cite{AlbinC}.} 
First note that, for any $k>0$
\BQNY
\pk{ Y(0)> u, Y(aq(u) k)> u}=\prod_{i=1}^n \pk{ X_i(0)> u, X_i(aq(u)k)> u}\le \prod_{i=1}^n \pk{X_i(0)+ X_i(aq(u)k)> 2 u}.
\EQNY
Therefore, 
similar arguments as in the proof of Lemma 1 therein imply
\BQNY
\limit{u} \frac{q(u)}{\pk{Y(0)> u}} \pk{ \sup_{k\in\{0,\ldots, \lfloor T/(aq(u))\rfloor\} }Y(aq(u)k) > u}&=& \frac{T}{a}
\pk{ \bigcap_{l=1}^\IF \{\mathcal{Z}(al) \le 0\}}
\EQNY
for any $a>0$. Moreover, the finiteness of the generalized Pickands constant $\H_{\alpha_1\ldot \alpha_n}(C_1\ldot C_n)$ and the asymptotic equation \eqref{eqMain} can be  established  as in \cite{AlbinC}, using the results of Lemmas \ref{AlbA}-\ref{AlbC}. In fact, $\H_{\alpha_1\ldot \alpha_n}(C_1\ldot C_n) >0$ follows directly from \eqref{lowerB}. This completes the proof. \QED

\prooftheo{Thm:Berman} The claim follows by checking the Assumptions 3.I and 3.II in Theorem 3.1 in \cite{Berman82}.
Assumption 3.I can be established with the aid of \nelem{AlbA},
where we have (with the notation \tb{as} in \cite{Berman82}) $w(u)=nu$, $v(u)=u^{2/\almin}$ and $Z(t)=\mathcal{Z}(t)$. Furthermore, it follows that
\BQNY
\lim_{d\to\IF}\limsup_{u\to\IF}v(u)\int_{d/v(u)}^t\pk{Y(s)>u|Y(0)>u}\ ds&\le& \lim_{d\to\IF}\limsup_{u\to\IF}v(u)\int_{d/v(u)}^t\pk{X_i(s)>u|X_i(0)>u}\ ds,
 \EQNY
 where $X_i$ is some of $X_i$'s such that $\alpha_i=\almin$. Therefore, Assumption 3.II
 can be verified as in Section 7 therein, and thus the proof is complete. \QED

\prooftheo{ThmB2} Initially we establish the proof for the case $j=1$. Introduce a new random process $Z$ defined by
$$
Z(t)=X_i(t-(i-1)T),\ \ \ t\in[(i-1)T,iT),\ \ 1\le i\le n.
$$
For any $u\ge0$ we have
\BQNY
\pk{\sup_{t \in [0,T]} X_{1:n}(t)> u}&=&\pk{\sup_{t \in [0,nT]} Z(t)> u}.
\EQNY
By the Bonferroni inequality and Pickands theorem (see Eq. \eqref{pic})
\BQN\label{eq:max1}
\pk{\sup_{t \in [0,nT]} Z(t)> u}&\le& \sum_{i=1}^n \pk{\sup_{t \in [0,T]} X_i(t)> u}\nonumber\\
&=&T n  \H_{\alpha} u^{\frac{2}{\alpha}} \Psi(u)(1+o(1)),\ \ u\to\IF
\EQN
and further
\BQN\label{eq:max2}
\pk{\sup_{t \in [0,nT]} Z(t)> u}&\ge& \sum_{i=1}^n \pk{\sup_{t \in [0,T]} X_i(t)> u}-\Sigma_1(u),
\EQN
with
$$\Sigma_1(u)=\sum_{1\le i<j\le n}\pk{\sup_{t \in [0,T]} X_i(t)> u}\pk{\sup_{t \in [0,T]} X_j(t)> u}.$$
Moreover, in view of \eqref{pic}
\BQN\label{eq:max3}
\Sigma_1(u)=o\left( u^{\frac{2}{\alpha}} \Psi(u)\right),\ \ u\to\IF.
\EQN
Consequently, the claim for the case $j=1$ follows from \eqref{eq:max1}-\eqref{eq:max3}. Next, we give only the proof of the case $j=n-1$ since the other cases follow by similar arguments. For notational simplicity denote
\BQNY
&&A_i(t,u)=\{X_1(t)>u,\ldots, X_{i-1}(t)>u,X_i(t)\le u, X_{i+1}(t)>u,\ldots, X_n(t)>u\},\ \ 1\le i \le n,\\
&&B(t,u)=\{X_1(t)>u,X_2(t)>u,\ldots, X_n(t)>u\}.
\EQNY
For any $u>0$ we have
\BQN\label{eq:Y1}
\pk{\sup_{t \in [0,T]} X_{n-1:n}(t)> u}&\le&\pk{\exists_{t \in [0, T]}\cup_{i=1}^nA_i(t,u)\cup B(t,u)}\nonumber\\
&\le& \pk{\exists_{t \in [0, T]}  B(t,u)} +\sum_{i=1}^n \pk{\exists_{t \in [0, T]} A_i(t,u) }\nonumber\\
&\le& \pk{\sup_{t \in [0,T]} \min_{1 \le i \le n} X_i(t)> u } +\sum_{i=1}^n \pk{\sup_{t \in [0,T]} \min_{1 \le j\le n, j\neq i} X_j(t)> u}.
\EQN
\gE{Further, for any $u>0$}
\BQN\label{eq:Y2}
\pk{\sup_{t \in [0,T]} X_{n-1:n}(t)> u}&\ge&\pk{\exists_{t \in [0, T]}\cup_{i=1}^nA_i(t,u) }\nonumber\\
&=&\pk{\cup_{i=1}^n\{\exists_{t \in [0, T]}A_i(t,u)\} }\nonumber\\
&\ge& \sum_{i=1}^n \pk{\sup_{t \in [0,T]} \min_{1 \le j\le n, j\neq i} X_j(t)> u}\pk{\sup_{t\in[0,T]}X_i(t)\le u} -\Sigma_2(u), 
\EQN
where
\BQNY
\Sigma_2(u)=\sum_{1\le i<j\le n}\pk{\exists_{t \in [0, T]} A_i(t,u), \exists_{s \in [0, T]} A_j(s,u) }.
\EQNY
By the independence of $X_i's$, we conclude that
\BQN\label{eq:Y3}
\Sigma_2(u)\le n^2 \prod_{i=1}^n \pk{\sup_{t\in[0,T]}X_i(t)> u}.
\EQN
Consequently, the claim follows from \eqref{eq:Y1}-\eqref{eq:Y3} and an application of \netheo{ThmB}. \QED

\cP{\prooftheo{ThmNS} Denote $\widetilde{Y}(t)=\min_{1\le i\le n}\widetilde{X_i}(t), t\ge0$, and let $w(u)=\sum_{i=1}^n b_i^2 u$.
As in the proof of \nelem{AlbA}
for any grid of points $0\le t_0 < t_1 < \cdots < t_d< \IF$
\BQNY
\Bigl( w(u)(\widetilde{Y}(q(u) t_1)- u) \ldot w(u)(\widetilde{Y}(q(u)t_d)- u) \Bigr) \Bigl \lvert (\widetilde{Y}(0)>u) &\todis & \sum_{i=1}^n b_i^2 \Bigl( \widetilde{\mathcal{Z}}(t_1) \ldot \widetilde{\mathcal{Z}}(t_d)\Bigr)
\EQNY
holds as $u\to\IF$. \gE{Results analogous to} \nelem{AlbB} and \nelem{AlbC} for $\widetilde{Y}$ can be derived with similar arguments as in the case of $Y$. Consequently, \gE{the proof is established by repeating the arguments in} the proof of \netheo{ThmB}. \QED  }

\bigskip
{\bf Acknowledgement}: We are thankful to the referee  for several suggestions which improved our manuscript. Support from Swiss National Science Foundation Project 200021-140633/1
and the project RARE -318984 (an FP7 Marie Curie IRSES Fellowship) is kindly acknowledged.
The first author also acknowledges partial support by NCN Grant No 2011/01/B/ST1/01521 (2011-2013).
E. Hashorva thanks Patrik Albin for kindly sending a copy of \cite{AlbinPHD}.

\bibliographystyle{plain}
\bibliography{supminCCCC}
\end{document}